\documentclass[10pt]{article}
\usepackage[utf8]{inputenc}
\usepackage{mathtools,
        mathrsfs,amssymb,amsthm, amsmath,
        bm}
\usepackage{hyperref}
\usepackage{url}
\usepackage[a4paper, left = 1in, right = 1in, top =1in, bottom = 1in]{geometry} 
\usepackage{dsfont}   
\usepackage[british]{babel}

\newtheorem{thm}{Theorem}
\newtheorem{cor}{Corollary}
\newtheorem{pro}{Proposition}

\def\[#1\]{\begin{align*}#1\end{align*}}

\newcommand{\ceq}{\coloneqq}

\begin{document}
\title{Fermat's Equation Has No Solution with Some Prime Components}

\author{Yu-Lin Chou}

\date{}

\maketitle
\begin{abstract}
Within the scope of elementary number theory, we prove that, as the
main result, if $1 \leq x < y < z$ are integers such that at least
one of $y, z, x+y$ is prime then $x^{n}+y^{n} \neq z^{n}$ for every
odd integer $n \geq 3$. This result covers a special case of a conjecture
of Abel, and furnishes a definite way to construct infinitely many
setwise coprime integers that do not satisfy the Fermat's equation
uniformly in $n$.\\

\end{abstract}

\section{Introduction}

Fermat's last theorem (FLT) may be paraphrased as the assertion that
for every odd integer\footnote{It is clear that one only needs to consider odd primes $\geq 3$. But our proof does not depend on this consideration.}
$n \geq 3$ the equation $x^{n}+y^{n}=z^{n}$ has no solution in setwise
coprime integers $x, y, z \geq 1$. 

In what follows, a triple of integers $x,y,z \geq 1$ is called a
\textit{Fermat's triple} iff $x,y,z$ are setwise coprime and $x^{n}+y^{n}\neq z^{n}$
for every odd integer $n \geq 3$.

Seeking after an elementary, hopefully short proof of FLT, we propose
some conditions on the indeterminates $x,y,z$ of the Fermat's equations
$x^{n}+y^{n}=z^{n}$ where $n \geq 3$ is an odd integer; specifically,
we give some sufficient conditions for a triple of setwise coprime
integers to be a Fermat's triple. 

We first prove an auxiliary result: \begin{pro} If $x,y,z \geq 1$ are setwise coprime integers, and if $(x+y)/\gcd(x+y,z)$ is coprime to $z$, then $(x,y,z)$ is a Fermat's triple. \qed \end{pro}\label{prop1} 

Abel conjectured \cite{r} that every triple of integers $\geq 1$
with some component being a prime power is a Fermat's triple; we prove
a result that covers a special case of the Abel's conjecture: \begin{thm} Every triple of setwise coprime integers $1 \leq x < y < z$ such that at least one of $y,z,x+y$ is prime is a Fermat's triple. \qed \end{thm}\label{thm1}
From Theorem \ref{thm1} we obtain a precise way to construct infinitely
many nontrivial triples of integers $x,y,z \geq 1$ such that $x^{n}+y^{n}\neq z^{n}$
for every odd integer $n \geq 3$; we give by applying Theorem \ref{thm1}
a constructive, short proof of the following result: \begin{cor} There are infinitely many Fermat's triples. \qed \end{cor}\label{cor1}

Some utterances regarding the ``methodology'' implied here: One
of the most famous, attractive elementary approaches to FLT would
be the sufficient condition, given by Sophie Germain at least 150
years ago, on the exponent $n$ for FLT to hold. Another intention
of the present paper is to ask how far one would travel, with elementary
considerations alone, by starting instead from properties of the Fermat's
(or non-Fermat's) triples that seem appealing both formally and empirically.

\section{Proofs}

\begin{proof}[Proof of Proposition 1]

We argue by contradiction; the main tool is the fact (under the assumption)
that $x+y$ divides $x^{n}+y^{n}$.

Indeed, from the hypothesis that $n \geq 3$ is an odd integer we
have the factorization\footnote{A quick justification may be obtained from the observation that $(x+y)S_{n} = \sum_{i=1}^{n}(-1)^{i-1}x^{i}y^{n-i} + \sum_{i=0}^{n-1}(-1)^{i}x^{i}y^{n-i}$.}
$x^{n}+y^{n} = (x+y)S_{n}$ where $S_{n} \ceq \sum_{i=0}^{n-1}(-1)^{i}x^{i}y^{n-1-i}$.
If $d \ceq \gcd(x+y, z)$, if $d_{1} \ceq z/d$, and if $d_{2} \ceq (x+y)/d$,
then \[
d_{1} z^{n-1} = \frac{z^{n}}{d} = \frac{x^{n}+y^{n}}{d} = d_{2}S_{n}.
\]Since $\gcd(d_{2}, z) = 1$ by assumption, we have $z^{n-1} \mid S_{n}$;
but $\gcd(d_{1},d_{2}) = 1$. \end{proof}

\begin{proof}[Proof of Theorem 1]

Suppose the statement is false, so that there is some counterexample
$(x,y,z)$. Let $d,d_{1},d_{2}$ be as in the proof of Proposition
\ref{prop1}.

Since $(x+y)^{n} > x^{n}+y^{n}=z^{n}$ by assumption, we have $x+y > z$.
If $x+y$ is prime, then, from the inequality $x+y > z$, it follows
that $d = 1$; so \[ \gcd\bigg( \frac{x+y}{d}, z \bigg) = \gcd(x+y,z) = d = 1.\]This
is by Proposition \ref{prop1} impossible.

Since $d=1$ implies that $(x+y)/d$ is coprime to $z$, we have $d \geq 2$
under the assumption. If $z$ is prime, then $d = z$. Moreover, since
$x+y < 2z$, which can be obtained by inspection or by an application
of H{\"o}lder's inequality\footnote{The proof does not depend on H{\"o}lder's inequality, although H{\"o}lder's inequality gives a slightly sharper upper bound, from $2$ to $2^{1-\frac{1}{n}}$.}
to the product $1\cdot x + 1 \cdot y$ with respect to counting measure,
we have \[ dd_{2} = zd_{2} = x+y < 2z.\]This implies that $d_{2} = 1$,
and hence \[ \gcd\bigg( \frac{x+y}{d}, z \bigg) = \gcd(d_{2},z) = \gcd(1, z) = 1; \]but
this is again impossible by Proposition \ref{prop1}. 

Suppose $y$ is prime. Since $y^{n}=z^{n}-x^{n}$ by assumption, we
have the apparent factorization $y^{n}=(z-x)S'_{n}$ for exactly one
integer $S'_{n} \geq 1$. Since $y > z-x$, if $z-x \geq 2$ then
$z-x$ does not divide $y^{n}$ and hence $S'_{n}$ cannot be an integer.
But $x < y < z$ under the assumption, so $z-x = 1$ implies $y = z$
or $=x$; it is thus impossible that $z-x = 1$. \end{proof}

\begin{proof}[Proof of Corollary 1]

The desired infinitude can be constructed in an elementary way: Given
any (sufficiently large, if the trivial cases are to be excluded in
the first place) prime $y \geq 2$, let $z > y$ be an integer coprime
to $y$. If $1 \leq x < y$ is an integer, then the integers $x,y,z$
are setwise coprime. Moreover, since $y$ is the second largest component
of the triple $(x,y,z)$, Theorem \ref{thm1} implies that $(x,y,z)$
is a Fermat's triple. Since there are infinitely many primes by the
Euclid's theorem, we are all set. \end{proof}

\end{document}